\documentclass[leqno,12pt]{article}
\usepackage{amssymb,amsfonts}
\usepackage{amsmath,latexsym}
\usepackage{french}
\usepackage{amstext,epic,eepic,epsf,pslatex}
\usepackage{graphicx}

\newcommand{\bepr}{{\em Proof} } 
\newcommand{\enpr}{\hfill \rule{.5em}{.5em}}

\newcommand{\R}{{\mathbb R}}

\newcommand{\Tr}{\hbox{Tr\,}}

\def\Xint#1{\mathchoice
{\XXint\displaystyle\textstyle{#1}}%
{\XXint\textstyle\scriptstyle{#1}}%
{\XXint\scriptstyle\scriptscriptstyle{#1}}%
{\XXint\scriptscriptstyle\scriptscriptstyle{#1}}%
\!\int}
\def\XXint#1#2#3{{\setbox0=\hbox{$#1{#2#3}{\int}$ }
\vcenter{\hbox{$#2#3$ }}\kern-.6\wd0}}

\newtheorem{prop}{Proposition}[section] 
\newtheorem{thm}{Theorem}[section]

\begin{document}

\title{Projective properties of Divergence-free symmetric tensors, and new dispersive estimates in gas dynamics}

\author{Denis Serre \\ \'Ecole Normale Sup\'erieure de Lyon\thanks{U.M.P.A., UMR CNRS--ENSL \# 5669. 46 all\'ee d'Italie, 69364 Lyon cedex 07. France. {\tt denis.serre@ens-lyon.fr}}}

\date{}

\maketitle

\begin{english}
\begin{abstract}
The class of Divergence-free symmetric tensors is ubiquitous in Continuum Mechanics. We show its invariance under projective transformations of the independent variables. This action, which preserves the positiveness, extends Sophus Lie's group analysis of Newtonian dynamics.

When applied to models of gas dynamics --~such as Euler system or Boltzmann equation,~-- in combination with Compensated Integrability, this yields new dispersive estimates. The most accurate one is obtained for mono-atomic gases. Then the space-time integral of $t\rho^\frac1d p$ is bounded in terms of the total mass and moment of inertia alone.
\end{abstract}
\end{english}



\paragraph{Notations.} The transpose of a rectangular matrix $Q$ is $Q^T$. If $u$ and $v$ are vectors, we denote $u\otimes v=uv^T$. The unit sphere of $\R^{1+d}$ is $S^d$. The Euclidian norm in $\R^n$ is $|\cdot|$. A {\em tensor} is a matrix-valued map $S(x)$ defined over some open domain of $\R^n$. If $S$ is $m\times n$ and its entries are distributions, one may take its divergence row-wise, which we denote ${\rm Div}\,S$ with a capital letter D~; it is an $m$-vector valued distribution. The tensor is {\em Divergence-free} if ${\rm Div}\,S=0$.

\section{Introduction}

In a series of articles \cite{Ser_DPT,Ser_JMPA,Ser_Prague,Ser_HS}, we studied the class of Divergence-free symmetric Tensors. These objects are ubiquitous in Continuum Mechanics, where they encode either the conservation of mass and momentum (classical mechanics) or energy and momentum (special relativity). The applications concern a vast list of models such as Euler or Navier-Stokes equations in compressible fluid dynamics, linear and nonlinear versions of the Maxwell's equations, kinetic (Boltzman) equations and mean-field (Vlasov-Poisson) models for plasmas or galaxies. They can also be employed for the study of scalar conservation laws \cite{DSLS}, or for systems of particles, such as hard spheres dynamics.

From a mathematical perspective, an interesting phenomenon happens whenever such a tensor is positive semi-definite, a property that results from the hypothesis that particles repel each other according to a radial force. For instance this applies to inviscid compressible gases, hard spheres dynamics or plasmas, but it doesn't to viscous gases (because the stress tensor is indefinite), to electro-magnetism (the magnetic component of the force is not radial) or to galaxies (gravity is attractive).
Positive semi-definiteness, plus the control of the Divergence in the space of (not necessarily vanishing) bounded measures imply a gain of integrability (see \cite{Ser_DPT,Ser_JMPA}). Say that the ambiant space is an open domain of $\R^n$, so that the tensor $S$ is $n\times n$. Then the expression $(\det S)^\frac1n$, which is naturally a locally bounded measure, is actually a measurable function of class $L^{\frac n{n-1}}$. This qualitative side, called {\em Compensated Integrability}, allows us to speak of the function  $(\det S)^\frac1{n-1}$. The quantitative side is an estimate of the latter in the form of a new functional inequality. For instance, if $S$ is compactly supported in $\R^n$, then
\begin{equation}
\label{eq:FI}
\left\| (\det S)^\frac1n\right\|_{\frac n{n-1}}\leqslant C_n\|{\rm Div}\,S\|_{\cal M},
\end{equation}
where the right-hand side involves the total mass of the measure $|{\rm Div}\,S|$. The absolute constant $C_n$ is sharp, as (\ref{eq:FI}) becomes an equality when $S(x)$ equals the identity matrix in a ball, and vanishes elsewhere. Inequality (\ref{eq:FI}) is a far-reaching extension of the isoperimetric inequality, as well as of the Gagliardo inequality. An alternate situation concerns Divergence-free positive tensors that are periodic with respect to a lattice~; if $\Gamma$ is a fundamental domain, then one has
\begin{equation}
\label{eq:nonJ}
\Xint-_\Gamma(\det S(z))^{\frac1{n-1}}dz\leqslant\left(\det\Xint-_\Gamma S(z)\,dz\right)^{\frac1{n-1}}.
\end{equation}
We emphasize that map $S\mapsto(\det S)^{\frac1{n-1}}$ is not concave, so that (\ref{eq:nonJ}) does not follow from Jensen's inequality for the non-concave, as it would do if the exponent $\frac1{n-1}$ was replaced by $\frac1n$ in both the left/right-hand sides.

\bigskip

The present paper is devoted to a new aspect of the theory. It has been known from the beginning that the class of Divergence-free symmetric tensors is stable under the action of the linear group through congruences,
$$S\longmapsto S_P(z):=PS(P^{-1}z)P^T,\qquad P\in {\bf GL}_n(\R).$$
We establish here that this class is actually stable under the action of the bigger group of projective transformations (Propositions \ref{p:projtr} and \ref{p:act}). Still, this action preserves the positiveness of tensors.

The applications to Continuum Mechanics involve the fact that physical models are described by a system of PDEs in the form of a Divergence-free symmetric tensor, together with one or several algebraic closure relations. In general the latter are not preserved under the projective action and therefore the model, whatever it be (Euler, Boltzman, etc), is not projectively invariant. It may happens however that special closure relations are preserved, so that some exceptional models enjoy this invariance, as observed already by Bobilev \& Ibragimov \cite{BI}, Bobylev \& Vilasi \cite{BV} and by Illner \cite{Illner}~:
\begin{itemize}
\item Isentropic Euler system ($\rho,u,p$ are the mass density, the velocity field and the pressure)
\begin{eqnarray*}
\partial_t\rho+{\rm div}_x(\rho u) & = & 0, \qquad (t,x)\in\R^{1+d}, \\
\partial_t(\rho u)+{\rm Div}_x(\rho u\otimes u)+\nabla_xp & = & 0
\end{eqnarray*}
for a mono-atomic gas:
$$p\equiv A\rho^{1+\frac2d}$$
where $A>0$ is a constant.
\item Non-isentropic Euler system, with the conservation law of energy ($e$ the specific internal energy)
$$\partial_t\left(\frac12\rho|u|^2+\rho e\right)+{\rm div}_x(\left(\frac12\rho|u|^2+\rho e+p\right)u)=0$$
and now
$$p\equiv\frac2d\,\rho e.$$
\item Vlasov equation
$$(\partial_t+v\cdot\nabla_x+F\cdot\nabla_v)f=0,\qquad (t,x,v)\in\R^{1+d+d}$$
in the unknown density $f(t,x,v)\ge0$, where the force $F(t,x)$ derives from a potential
$$F=-\nabla_xU[\rho],\qquad\rho(t,x):=\int_{\R^d}f(t,x,v)\,dv$$
and the particles interact through the Calogero--Moser potential $A|x|^{-2}$,
$$U[\rho](t,x)=A\int_{\R^d}\frac{\rho(t,y)}{|x-y|^2}\,dy.$$
\item Boltzman equation (see Section \ref{s:Bolt}) when, again, the particles interact through the Calogero--Moser potential.
\end{itemize} 
When it applies, this additional invariance is always associated with an extra conservation law, which involves the moment of inertia.

Whether a model is projectively invariant or not, the use of projective transformations, combined with Compensated Integrability, provides us with a one-parameter family of estimates, much richer than the single one established in \cite{Ser_DPT}. Optimizing the choice of the parameter, we obtain a new estimate, whose interest is two-fold. On the one hand, it is significantly sharper for large times, giving a better dispersion of the mass as $t\rightarrow+\infty$. On the other hand, it highlights the symmetry played by the total mechanical energy and the moment of inertia. This confirms that the assumptions about the initial data, in the theory of renormalized solutions for the Boltzman equation, are natural. 

\bigskip

Let us illustrate all this by a single example, taken from the dynamics of a mono-atomic gas in space dimension $d$ (Theorem \ref{th:inertia})~: the density $\rho$ and the pressure $p$ are estimated by
$$\int_0^{t_{\max}}t\,dt\int_{\R^d}\rho^\frac1dp\,dx\leqslant c_dM^\frac1d\left(\frac14\int_{\R^d}\!\int_{\R^d}\rho_0(x)\rho_0(x')|x'-x|^2dx\,dx'\right)^\frac12\,,$$
where $M$ is the total mass and $\rho_0$ the mass density at time $t=0$. The constant $c_d$ does not depend upon the solution or the time interval: when the flow is global-in-time, the estimate is valid with $t_{\max}=+\infty$. Amazingly, this estimate involves only the distribution of mass at initial time, but not the prescribed initial velocity field and temperature, in spite of the coupling during the evolution~! This sheds another light upon the {\em source solutions} of the Euler system, those for which $\rho_0$ is a Dirac mass.

\paragraph{Plan of the paper.} The abstract analysis is done in Section \ref{s:abst}. A first application, to either isentropic or non-isentropic gas dynamics, is given in Section \ref{s:gas}.
After a short account of the theory of Compensated Integrability, we carry the calculations for ideal gases, which culminate with Theorem \ref{th:inertia}. Section \ref{s:Bolt} is dedicated to the Boltzman equation, for which similar results hold true.

\paragraph{Warning.} The results stated  in Sections \ref{s:gas} and \ref{s:Bolt} concern flows in the entire space $\R^d$. Our admissibility criterion is that not only the total mass $M$, mechanical energy $E(t)$ and moment of inertia $I(t)$ are finite, but $M$ is constant, $t\mapsto E(t)$ is non-increasing, and the differential inequality (\ref{eq:tren}) is valid. In particular, we exclude the wild solutions constructed by De Lellis \& Sz\'ekelyhidi \cite{DeLSz}.

\paragraph{Acknowledgement.} I am indebted to Reinhard Illner (Univ. of Victoria, CA) for driving my attention towards important references.

\section{Projective transformations and Divergence-free symmetric tensors}\label{s:abst}

\subsection{State of the art}

A conservation law in $d+1$ space-time dimension,
\begin{equation}
\label{eq:single}
\partial_t\rho+{\rm div}_xm=0,
\end{equation}
expresses that some differential form $\omega$ of degree $d$ is closed. Performing an smooth change of variables $(s,y)=\phi(t,x)$, we rewrite (\ref{eq:single}) as another conservation law
\begin{equation}
\label{eq:sbar}
\partial_s\bar\rho+{\rm div}_y\bar m=0,
\end{equation}
where $(\bar\rho,\bar m)$ are obtained from $(\rho,m)$ through a linear transformation with variable coefficients. Distributional solutions of (\ref{eq:single}) yield distributional solutions of (\ref{eq:sbar}) because they both express the same property ${\rm d}\omega=0$, though in different coordinates.

Let us give an example, which turns out to be fundamental in the applications. The space-time domain being $(0,t_{\max})\times\R^d$, we choose a projective transformation
\begin{equation}
\label{eq:sytx}
s=\frac t{1+t\alpha}\,,\qquad y=\frac x{1+t\alpha}
\end{equation}
where $\alpha>0$ is some constant parameter. Then the dependent variables are transformed according to
$$\bar\rho:=(1+t\alpha)^d\rho,\qquad\bar m=(1+t\alpha)^{d+1}m-\alpha(1+t\alpha)^d\rho x,$$
where we point out that $1+t\alpha=(1-s\alpha)^{-1}$. 
When $\rho$ is positive, and thus plays the role of a mass density, it is meaningful to introduce a `velocity' field by $u:=\frac m\rho$\,. Then the new `momentum' $\bar m$ is given as $\bar\rho v$, where the new velocity is defined by $v=(1+t\alpha)u-\alpha x$.

\bigskip

Projective transformations such as (\ref{eq:sytx}) are meaningful in Classical dynamics too. They leave the simplest ODE
$$\frac{d^2x}{dt^2}=0$$
invariant, just because they transform lines of the $(t,x)$-space into lines of the $(s,y)$-space. Amazingly enough, the nonlinear ODE
$$\frac{d^2x}{dt^2}=-\nabla_xU(|x|),\qquad U(r):=\frac{\rm cst}{r^2}\,,$$
where $U$ is the so-called {\em Calogero--Moser potential}, is also invariant under the action of (\ref{eq:sytx}).

\bigskip

A.V. Bobylev, with either N. Kh. Ibragimov \cite{BI} or G. Valesi \cite{BV}, studied in a systematic way the action of projective transformations upon PDEs from Mathematical Physics. They observed that the Euler system of a mono-atomic gas is left invariant (see also \cite{Ser_AIF}), as well as the Boltzman equation when the particles interact through the Calogero--Moser potential. R. Illner \cite{Illner} proved the same properties for the Vlasov equation with the C.-M. potential. In all cases, this additional symmetry is associated with an extra conservation law, a property which corroborates N{\oe}ther's Theorem.

\subsection{A universal structure in Continuum Mechanics}

Several models of Continuum Mechanics can be decomposed into two parts. On the one hand, they share a couple of fundamental conservation laws, that of mass and linear momentum. These are expressed as a linear system of PDEs
\begin{equation}
\label{eq:DivS}
{\rm Div}_{t,x}S=0,
\end{equation}
where $S$ is a symmetric tensor whose entries are distributions, and the divergence is taken row-wise:
$$\forall i=0,\ldots,d,\qquad({\rm Div}_{t,x}\,S)_i:=\partial_ts_{i0}+\sum_{j=1}^d\partial_{x_j}s_{ij}=0.$$
Typically
$$S=\begin{pmatrix} \rho & m^T \\ m & \frac{m\otimes m}\rho+\sigma \end{pmatrix}$$
where $\rho, m,\sigma$ are the mass density, the linear momentum and the stress tensor. Several examples were described in \cite{Ser_DPT,Ser_JMPA}. On the other hand, each model is closed by algebraic, or differential-algebraic closure relations. Remark that from this point of view, the conservation of energy (when it applies) stands as a closure relation.

We wish to focus onto the somehow universal governing equations (\ref{eq:DivS}). It turns out that this structure is projectively invariant.

\begin{quote}
\begin{prop}\label{p:projtr}
Let $S(t,x)$ be a Divergence-free symmetric tensor. Let us write it blockwise
$$S=\begin{pmatrix} \rho & m^T \\ m & T \end{pmatrix},$$
where $\rho$ is scalar and thus $T$ is $d\times d$ and symmetric.

Let $\alpha$ be a positive constant, and consider the change of variables defined by (\ref{eq:sytx}). Then the symmetric tensor $\bar S$, whose blocks are defined by
\begin{eqnarray*}
\bar\rho  & = &  (1+t\alpha)^d\rho, \\
\bar m & = & (1+t\alpha)^{d+1}m-\alpha(1+t\alpha)^d\rho x , \\
\bar T & = & (1+t\alpha)^{d+2}T-\alpha(1+t\alpha)^{d+1}(m\otimes x+x\otimes m)+\alpha^2(1+t\alpha)^d\rho x\otimes x ,
\end{eqnarray*}
is Divergence-free in the $(s,y)$ variables defined by (\ref{eq:sytx}):
$${\rm Div}_{s,y}\bar S=0.$$
\end{prop}
\end{quote}

\paragraph{Remarks.} 1) The first row is transformed exactly as in our very first example. 2) The proof of the Proposition is elementary and follows from the formul\ae
$$\partial_t=(1-s\alpha)^2\partial_s-\alpha(1-s\alpha)y\cdot\nabla_y,\qquad\nabla_x=(1-s\alpha)\nabla_y.$$
3) The positivity of symmetric matrices is preserved by the transformations above. We have
$$\bar S=(1+t\alpha)^dPSP^T,\qquad P:=\begin{pmatrix} 1 & 0 \\ -\alpha x & (1+t\alpha)I_d \end{pmatrix}.$$

\bigskip

A more general result concerns the case where the divergence is non-zero.

\begin{quote}
\begin{prop}
\label{p:divcontrol}
With the notations above, but allowing ${\rm Div}_{t,x}S$ to be a non-zero distribution, we have
\begin{eqnarray*}
\langle\partial_s\bar\rho+{\rm div}_y\bar m,\phi\rangle & = & \langle\partial_t\rho+{\rm div}_x m,\psi\rangle , \\
\langle\partial_s\bar m+{\rm Div}_y\bar T,\vec A\rangle & = & \langle\partial_t m+{\rm Div}_xT,(1+\alpha t)\vec B\rangle -\alpha\langle\partial_t\rho+{\rm div}_xm,x\cdot\vec B\rangle 
\end{eqnarray*}
where the test functions are related by $\psi(t,x)=\phi(s,y)$ and $\vec B(t,x)=\vec A(s,y)$.

In particular, if ${\rm Div}_{t,x}S$ is a vector-valued bounded measure, then so is ${\rm Div}_{s,y}\bar S$ and we have
\begin{eqnarray*}
\|\partial_s\bar\rho+{\rm div}_y\bar m\|_{\cal M} & = & \|\partial_t\rho+{\rm div}_x m\|_{\cal M} , \\
\|\partial_s\bar m+{\rm Div}_y\bar T\|_{\cal M} & \leqslant & \|(1+\alpha t)(\partial_t m+{\rm Div}_xT)\|_{\cal M} +\alpha\|\,|x|(\partial_t\rho+{\rm div}_x m)\|_{\cal M} .
\end{eqnarray*}
\end{prop}
\end{quote}

The proof of the differential identities are straightforward and left to the reader. The estimates of the masses of measures are obtained by taking the supremum over those $\phi$ or $\vec A$ such that $|\phi(s,y)|\le1$ of $|\vec A(s,y)|\le1$ pointwise.

\subsection{The action of the projective group}

At first glance, the factors $(1+t\alpha)^d$ to $(1+t\alpha)^{d+2}$ in the definition of $\bar S$ might look weird. Besides, one could ask oneself why the Proposition holds true for the projective transformation (\ref{eq:sytx}), but does not for a general change of variable. The explanation of both facts comes from the following observation.

\begin{quote}
\begin{prop}\label{p:restr}
Let $\Sigma(\lambda,z)$ be a symmetric Divergence-free tensor over an open cone $\Gamma\subset\R^{1+n}$. We assume that $\Sigma$ is positively homogeneous of degree $-n-1$. Let us write it blockwise
$$\Sigma(\lambda,z)=\lambda^{-n-1}\begin{pmatrix} h\left(\frac z\lambda\right) & Z\left(\frac z\lambda\right) ^T \\ Z\left(\frac z\lambda\right) & H\left(\frac z\lambda\right) \end{pmatrix}.$$
Then the $n\times n$ symmetric tensor
$$\widetilde H(x)=H(x)-x\otimes Z(x)-Z(x)\otimes x+h(x)x\otimes x$$
is Divergence-free: ${\rm Div}_x\widetilde H=0$.
\end{prop}
\end{quote}

\bigskip

\bepr

When expressing ${\rm Div}_{\lambda,z}\Sigma=0$, we obtain
$${\rm div}_xZ=(x\cdot\nabla_x)h+(n+1)h,\qquad {\rm Div}_xH=(x\cdot\nabla_x)Z+(n+1)Z.$$
On the other hand, we always have
\begin{eqnarray*}
{\rm Div}_x(xZ^T+Zx^T) & = & (x\cdot\nabla_x)Z+({\rm div}_xZ)x+(n+1)Z, \\
{\rm Div}_x(hxx^T) & = & (x\cdot\nabla_xh)x+(n+1)hx.
\end{eqnarray*}
Making a linear combination of the four identities, we eliminate the zero-order terms and obtain the desired conclusion.

\enpr

\bigskip

Once again, the positivity of symmetric matrices is preserved by the transformation above:
$$\mu^T\widetilde H\mu=(-x\cdot\mu\,,\,\mu^T)\Sigma\binom{-x\cdot\mu}\mu,\qquad\forall\,\mu\in\R^n.$$

\bigskip

The strategy to pass from Proposition \ref{p:restr} to Proposition \ref{p:projtr} is to set $n=d+1$ and to extend $S$ as an $(n+1)$-dimensional tensor by
$$\Xi(\lambda,z)=\lambda^{-n-1}\begin{pmatrix} 0 & 0 \\ 0 & S\left(\frac z\lambda\right) \end{pmatrix}.$$
The tensor $\Xi$ is Divergence-free (obvious), homogeneous of degree $-n-1$. When $P\in {\bf GL}_{n+1}(\R)$, a congruence
$$\Sigma(w)=(\det P)^{-1}|\det P|^{-\frac2{n+1}}P\Xi(P^{-1}w)P^T$$
defines another symmetric Divergence-free tensor (Lemma 1.1 of \cite{Ser_DPT}), still homogeneous. One can then apply Proposition \ref{p:restr} to $\Sigma$. An appropriate choice of $P$ yields the transformation treated in Proposition \ref{p:projtr}.

\bigskip

The procedure $S\mapsto\Xi\mapsto\Sigma\mapsto\widetilde H$ defines an action of the linear group ${\bf GL}_{n+1}(\R)$ over the space of $n\times n$ Divergence-free symmetric tensors. Observing that a homothety $P=aI_{n+1}$ yields $\Sigma=\Xi$ and thus $\widetilde H=S$, we may pass to the quotient and state

\begin{quote}
\begin{prop}\label{p:act}
The composition of the linear maps $S\mapsto\Xi$, $\Xi\mapsto\Sigma$ (this one being an action of ${\bf GL}_{n+1}(\R)$) and $\Sigma\mapsto\widetilde H$ defines an action of the projective linear group ${\bf PGL}(n+1;\R)$ over the space of $n\times n$ Divergence-free symmetric tensors. 

This action preserves the positive semi-definiteness.
\end{prop}
\end{quote}

\subsection{Special Divergence-free tensors}

We recall (see \cite{Ser_DPT}) that among the class of positive Divergence-free tensors (DPTs), one encounters the (non-linear) class of {\em special DPTs}, which are cofactor matrices of Hessians with convex potentials:
\begin{equation}
\label{eq:special}
S=\widehat{{\rm D}^2_{t,x}\theta}.
\end{equation}

\begin{quote}
\begin{prop}
\label{p:special}
Let $S$ be as in (\ref{eq:special}), where $\theta$ is a convex function. Let $\alpha>0$ be a parameter. Then the Divergence-free tensor $\bar S$ defined in Proposition \ref{p:projtr} is again a special DPT,
$$\bar S=\widehat{{\rm D}^2_{s,y}\bar\theta},$$
where the new potential, given by
$$\bar\theta(s,y)=(1-\alpha s)\theta\left(\frac s{1-\alpha s}\,,\frac y{1-\alpha s}\right),$$
is convex.
\end{prop}
\end{quote}

\bigskip

\bepr

On the one hand, the function $\bar\theta$ is convex because $\theta$ can be written as the supremum of a collection of affine functions $(t,x)\mapsto pt+\xi\cdot x+c$ and thus $\bar\theta$ is the supremum of the affine functions $(s,y)\mapsto ps+\xi\cdot y+c(1-\alpha s)$.

On the other hand, a cumbersome though elementary calculation gives
$${\rm D}^2_{s,y}\bar\theta=(1+\alpha t)Q^T({\rm D}^2_{t,x}\theta) Q,\qquad Q:=\begin{pmatrix} 1+\alpha t & 0 \\ \alpha y  & I_d \end{pmatrix}.$$
One concludes by using the formula $\widehat{AB}=\widehat A \widehat B$, plus the identity
$$\widehat Q=\begin{pmatrix} 1 & -\alpha y^T \\ 0 & (1+\alpha t)I_d \end{pmatrix}.$$

\enpr

\subsection{Determinantal masses}

Let us recall (see \cite{Ser_JMPA}) the rigidity result, that if a Divergence-free positive tensor is homogeneous of degree $-d$ about a point $p=(t_0,x_0)\in\R^{1+d}$, then it is of the form
$$S=\mu\left(\frac z{|z|}\right)\,\frac{z\otimes z}{|z|^{d+2}}$$
where $z=\binom tx-p$ and $\mu$ is a non-negative measure over the unit sphere $S^d$. This was completed in \cite{Ser_HS} with the observation that, thanks to Pogorelov's Theorem (see \cite{Pog}), $S$ is actually
special in the sense of (\ref{eq:special}), and its potential $\theta$ is positively homogeneous of degree $1$. The converse is obviously true: if $\theta$ is positively homogeneous of degree $1$, then $\widehat{{\rm D}^2\theta}$ is Divergence-free and positively homogeneous of degree $-d$.

We showed in \cite{Ser_HS} that in this situation, the expression $(\det S)^\frac1d$ concentrates as a Dirac mass ${\rm Dm}(S;p)\delta_{(t,x)=p}$. We identified its weight, called a {\em Determinantal Mass}, as the volume of the convex body surrounded by the image\footnote{The gradient is positively homogeneous of degree zero. Its image is a closed hypersurface.} of $\nabla_{t,x}\theta$. This generalizes the observation that for smooth potentials, 
$$\left(\det\widehat{{\rm D}^2\theta}\right)^\frac1d=\det{\rm D}^2\theta$$
is the Jacobian of the map $(t,x)\mapsto\nabla_{t,x}\theta$.

\bigskip

This homogeneity of degree one (for the potential) is actually preserved by our projective transformation: Let us start from the Euler identity
$$(t-t_0)\partial_t\theta+(x-x_0)\cdot\nabla_x\theta=\theta.$$
Combining with the expressions
$$\partial_s\bar\theta=-\alpha\theta+\frac1{1-\alpha s}\,\partial_t\theta+\frac{\alpha y}{1-\alpha s}\,\cdot\nabla_x\theta,\qquad\nabla_y\bar\theta=\nabla_x\theta,$$
we have first 
$$\partial_s\bar\theta=(1+\alpha t_0)\partial_t\theta+\alpha x_0\cdot\nabla_x\theta$$
and we derive
$$(s-s_0)\partial_s\bar\theta+(y-y_0)\cdot\nabla_y\bar\theta=\bar\theta,$$
where $(s_0,y_0)$ is the corresponding point in the $(s,y)$ coordinates. 

Now, because
$\nabla_{s,y}\bar\theta(s,y)$ is the composition of $\nabla_{t,x}\theta(t,x)$ with the linear map of matrix
$$A=\begin{pmatrix} 1+\alpha t_0 & \alpha x_0^T \\ 0 & I_d \end{pmatrix},$$
the volumes surrounded by the image of $\nabla_{s,y}\bar\theta$ or by that of $\nabla_{t,x}\theta$ differ from each other by the factor $\det A=1+\alpha t_0$. Let us summarize this analysis:

\begin{quote}
\begin{thm}\label{th:hom}
Let the symmetric, positive Divergence-free tensor $S$ be homogeneous of degree $-d$ about a point $p=(t_0,x_0)$. Then under the projective change of variables, the tensor $\bar S$ has the same homogeneity about $\bar p$, the image of $p$, and the determinantal masses are related through
\begin{equation}
\label{eq:DMs}
{\rm Dm}(\bar S;\bar p)=(1+\alpha t_0){\rm Dm}(S;p).
\end{equation}
\end{thm}
\end{quote}

\section{Application to gas dynamics}\label{s:gas}

When considering the evolution of a fluid, $d$ is the space dimension and $S$ is the mass-momentum tensor, with $\rho,m$ being the mass density and linear momentum. The tensor $\sigma=T-\frac{mm^T}\rho$\,, which is the Schur complement of $\rho$ in $S$, is the (opposite of the) stress tensor. In the transformation $S\mapsto \bar S$, we have
$$\bar\sigma=(1+t\alpha)^{d+2}\sigma.$$
Remark that for an inviscid gas, where $\sigma=pI_d$ ($p\ge0$ the pressure), $S$ is positive semi-definite (definite whenever $\rho>0$ and $p>0$), with determinant $\rho p^d$.

In general, there is no reason why $\bar S$ should be the mass-momentum tensor of another gas flow. For a barotropic gas, this happens only if
$$(1+t\alpha)^{d+2}p(\rho)\equiv p\left((1+t\alpha)^d\rho\right).$$
This means that $p$ is homogeneous of degree 
$$\gamma=\gamma_d:=1+\frac2d\,,$$
which corresponds to the so-called {\em mono-atomic gas}. The symmetry group of the PDEs has then an extra dimension.

\begin{quote}
\begin{prop}[\cite{BI,BV,Ser_AIF}]\label{p:xtramono}
The Euler system for a mono-atomic inviscid barotropic gas is invariant under the transformation 
$$(x,t,\rho,u)\mapsto\left(\frac x{1+t\alpha}\,,\frac t{1+t\alpha}\,,(1+t\alpha)^d\rho,(1+t\alpha)u-\alpha x\right).$$
\end{prop}
\end{quote}

\subsection{Full gas dynamics}

In full gas dynamics, the flow is described by the triplet $(\rho,u,e)$ where $(\rho,u)$ are as above, and $e$ is the specific internal energy. It is governed by $d+2$ conservation laws, namely those of mass, momentum and energy. The two first write as ${\rm Div}_{t,x}S=0$ where $S$ is the mass-momentum tensor, while the last one is
\begin{equation}
\label{eq:consen}
\partial_t\left(\frac12\rho|u|^2+\rho e\right)+{\rm div}_x\left(\left(\frac12\rho|u|^2+\rho e+p\right)u\right)=0.
\end{equation}
Under mild integrability assumptions, the total mass, momentum and energy
$$M=\int_{\R^d}\rho(t,x)\,dx,\qquad Q=\int_{\R^d}(\rho u)(t,x)\,dx,\qquad E=\int_{\R^d}\left(\frac12\rho|u|^2+\rho e\right)(t,x)\,dx$$
are constants of the motion. 

We assume the equation of state of a perfect gas:
$$p=(\gamma-1)\rho e,$$
where $\gamma>1$  still denotes the adiabatic constant.
Defining
$$\bar e=(1+t\alpha)^2e,\qquad \bar p=(1+t\alpha)^{d+2}p$$
the effect of the projective transformation over (\ref{eq:consen}) is that
\begin{equation}
\label{eq:tren}
\partial_s\left(\frac12\bar\rho|v|^2+\bar\rho\bar e\right)+{\rm div}_x\left(\left(\frac12\bar\rho|v|^2+\bar\rho\bar e+\bar p\right)v\right)=\frac{d \alpha}{1-s\alpha}\, (\gamma_d-\gamma)\rho e.
\end{equation}
This time, the new field does satisfy the equation of state: $\bar p=(\gamma-1)\bar\rho\,\bar e$, but this happens to the expense of a non-trivial right-hand side in (\ref{eq:tren})~: the field $(\bar\rho,v,\bar e)$ does not satisfy the conservation law of energy.
Once again, the only case where the right-hand side vanishes is for $\gamma=\gamma_d$.

\begin{quote}
\begin{prop}\label{p:xtraonom}
The Euler system for a mono-atomic inviscid gas is invariant under the transformation 
$$(x, t,\rho,u,e)\mapsto\left(\frac x{1+t\alpha}\,,\frac t{1+t\alpha}\,,(1+t\alpha)^d\rho,(1+t\alpha)u-\alpha x,(1+t\alpha)^2e\right).$$
\end{prop}
\end{quote}

\subsection{Compensated Integrability}

Let us recall Theorem 2.3 of \cite{Ser_Prague}:

\begin{quote}
\begin{thm}\label{th:prague}
Let $S$ be a positive semi-definite Divergence-free symmetric tensor over $Q_\tau:=(0,\tau)\times\R^d$, written blockwise
$$S=\begin{pmatrix} \rho & m^T \\ m & T \end{pmatrix}.$$
We assume that $S$ is integrable over $Q_\tau$, and that its (well-defined) normal traces $S\vec e_t$ at the top ($t=\tau$) and bottom ($t=0$) are bounded measures (for instance, integrable over $\R^d$).

Then $(\det S)^{\frac1{d+1}}$, which is {\em a priori} a bounded measure, actually belongs to $L^{1+\frac1d}(Q_\tau)$, and we have
\begin{equation}
\label{eq:precis}
\int_0^\tau dt\int_{\R^d}(\det S)^\frac1ddx\leqslant\frac1d\,\left(\frac{2(d+1)}{|S^d|}\right)^\frac1dM^\frac1d(\|m(0,\cdot)\|_{\cal M}+\|m(\tau,\cdot)\|_{\cal M}),
\end{equation}
with 
$$M:=\int_{\R^d}\rho(t,x)\,dx$$
the constant total mass, and $|S^d|$ is the area of the unit sphere of $\R^{d+1}$.
\end{thm}
\end{quote}

\bigskip

The constant appearing in the right-hand side of (\ref{eq:precis}) is almost sharp, but its value is not important for us. In the sequel, upper bounds will involve other constants, depending only upon the space dimension, which will always be denoted $c_d$. We emphasize that these constants are {\bf moderate}. Typically, they satisfy $c_3<10$.

\bigskip

When Theorem \ref{th:prague} is applied to a gas, we have $m=\rho u$ and $T=\rho u\otimes u+p I_d$, so that $\det S=\rho p^d$. Cauchy--Schwarz inequality yields
$$\|m(t,\cdot)\|_{\cal M}^2\leqslant M\int_{\R^d}\rho|u|^2dx=2M E_{kin}(t)\le2ME(t)$$
where $E_{kin}(t),E(t)$ denote the kinetic and total mechanical energies at time $t$. Whether the gas is isentropic or not, the total energy of an admissible flow is non-increasing or constant. We infer the estimate
\begin{equation}
\label{eq:estihp}
\int_0^\tau dt\int_{\R^d}\rho^\frac1dp\,dx\leqslant c_dM^{\frac1d}\sqrt{ME(0)\,}\,,
\end{equation}
first established in \cite{Ser_DPT}. We point out that the right-hand side does not depend upon the length $\tau$ of the time interval. Therefore (\ref{eq:estihp}) holds true even when the flow is globally defined:
$$\int_0^\infty dt\int_{\R^d}\rho^\frac1dp\,dx\leqslant c_dM^{\frac1d}\sqrt{ME(0)\,}\,.$$

\subsection{Mono-atomic gas}\label{p:mono}

Let $(\rho,u,e)$ be an admissible flow of a mono-atomic gas, defined over $(0,t_{\max})\times\R^d$, with finite mass and energy. It may be either isentropic or not. 

Given  a constant $\alpha>0$, we perform the transformation studied above, in either Proposition \ref{p:xtramono} or \ref{p:xtraonom}. What matters here is that the transformed field is still a flow, and an admissible one, because of the mono-atomic assumption.
This new flow is defined over $Q_\tau$, where 
$$\tau=\tau_\alpha:=\frac{t_{\max}}{1+t_{\max}\alpha}\,.$$
It has the same mass as the original one:
$$\int_{\R^d}\bar\rho(s,y)\,dy=\int_{\R^d}(1+t\alpha)^d\rho(t,(1+t\alpha)y)\,dy=\int_{\R^d}\rho(t,x)\,dx\equiv M.$$
Assuming in addition that $\rho(0,\cdot)\in L^1(|x|^2dx)$ -- the moment of inertia is finite, -- then its total energy
$$E_\alpha(s)=\int_{\R^d}\left(\frac12\,\bar\rho\,|v|^2+\bar\rho\,\bar e\right)dy=\int_{\R^d}\left(\frac12\,\rho|(1+t\alpha)u-\alpha x|^2+(1+t\alpha)^2\rho e\right)dx$$
is finite too. Applying (\ref{eq:estihp}) to the new flow, we obtain
\begin{equation}
\label{eq:baralp}
\int_0^{\tau_\alpha}ds\int_{\R^d}\bar\rho^\frac1d\bar p\,dy\leqslant c_dM^\frac1d\sqrt{ME_\alpha(0)\,}\,.
\end{equation}
We calculate on the one hand
\begin{eqnarray*}
\int_0^{\tau_\alpha}ds\int_{\R^d}\bar\rho^\frac1d\bar p\,dy & =  & \int_0^{t_{\max}}\frac{dt}{(1+t\alpha)^2}\int_{\R^d}(1+t\alpha)^{d+3}\rho^\frac1dp\,\frac{dx}{(1+t\alpha)^d} \\
\nonumber
& = & \int_0^{t_{\max}}(1+t\alpha)dt\int_{\R^d}\rho^\frac1dp\,dx.
\end{eqnarray*}
On the other hand, we have 
$$E_\alpha(0)=\int_{\R^d}\left(\frac12\,\rho(0)|u(0)-\alpha x|^2+\rho(0) e(0)\right)\,dx.$$
Inserting these into (\ref{eq:baralp}), we infer
$$\int_0^{t_{\max}}t\,dt\int_{\R^d}\rho^\frac1dp\,dx \leqslant   c_dM^\frac1d\left(\inf_{\alpha>0}\,\frac M{\alpha^2}\,E_\alpha(0)\right)^\frac12.$$ 
We can therefore conclude
\begin{equation}
\label{eq:noninv}
\int_0^{t_{\max}}t\,dt\int_{\R^d}\rho^\frac1dp\,dx\leqslant  c_dM^\frac1d\left(\lim_{\alpha\rightarrow+\infty}\,\frac M{\alpha^2}\,E_\alpha(0)\right)^\frac12 =  c_dM^\frac1d\sqrt{MI(0)\,}\,,
\end{equation}
where
$$I(0)=\int_{\R^d}\rho(0,x)\,\frac{|x|^2}2\,dx$$
is the moment of inertia at initial time. 

Remark that the left-hand side of (\ref{eq:noninv}) is translation invariant, while its right-hand side is not. Choosing a different origin of the physical space, we may replace the quantity $MI_(0)$ by
$$M\int_{\R^d}\rho(0,x)\,\frac{|x-\hat x|^2}2\,dx,$$
where $\hat x\in\R^d$ is constant.
Taking the infimum as $\hat x$ runs over $\R^d$, we obtain the following statement.

\begin{quote}
\begin{thm}\label{th:inertia}
Consider the admissible flow of a mono-atomic gas (either isentropic or not) in $(0,t_{\max})\times\R^d$ with finite mass, energy and moment of inertia at time $t=0$. Then we have
\begin{equation}
\label{eq:inertia}
\int_0^{t_{\max}}t\,dt\int_{\R^d}\rho^\frac1dp\,dx\leqslant  c_dM^\frac1d\left(\frac14\int_{\R^d}\!\int_{\R^d}\rho(0,x)\rho(0,x')|x'-x|^2dx\,dx'\right)^\frac12.
\end{equation}
for some absolute constant $c_d$, where $\rho(0,\cdot)$ is the mass density at inital time.
\end{thm}
\end{quote}

\paragraph{Comments.}
\begin{enumerate}
\item At first glance, Estimate (\ref{eq:inertia}) suggests that the hypothesis of finite energy might be useless, since it does not involve the initial velocity field. This is a deadly false impression: Compensated Integrability is valid only if the mass-momentum tensor is integrable over $(0,t_{\max})\times\R^d$, or at least over $(\epsilon,t^*)\times\R^d$ for every $0<\epsilon<t^*<t_{\max}$. Because of the positivity, this amounts to saying that its trace is integrable, which means that 
$$\int_\epsilon^{t^*}(M+E(t))dt<\infty.$$
In the non-isentropic case, this is equivalent to assuming $E(0)<\infty$, while in the isentropic case, where $t\mapsto E(t)$ is non-increasing, it says that the total energy is finite whenever $t>0$. To avoid complications, we choose to assume a finite initial energy.
\item Both sides of (\ref{eq:inertia}) have the same physical dimension $M^{1+\frac1d}L$.
\item We already new from (\ref{eq:estihp}) that the function
$$\pi(t):=\int_{\R^d}\rho^\frac1dp\,dx$$
is integrable over $(0,t_{\max})$. The main information carried by the new estimate is that, if $t_{\max}=+\infty$, then $t\pi(t)$ is integrable at infinity. This narrows the gap between Estimate (\ref{eq:estihp}) and the decay observed for the global classical solutions constructed in \cite{Ser_AIF}.
\item As a corollary, we get the well-known fact that `source-solutions', for which the whole mass concentrates in one point at initial time, must have an infinite energy. A more elementary proof  follows from the quite obvious functional inequality
\begin{equation}
\label{eq:uncert}
\left(\int_{\R^d}g(x)\,dx\right)^{3+\frac2d}\leqslant c_d\int_{\R^d}\!\int_{\R^d}g(x)g(x')|x'-x|^2dx\,dx'\cdot\int_{\R^d}g(x)^{1+\frac2d}dx
\end{equation}
for non-negative functions (see Appendix), applied to $\rho(t,\cdot)$. It shows that at constant mass, if either the internal energy or the moment of inertia tends to zero, then the other one tends to infinity --~a kind of uncertainty principle.
\item
Estimate (\ref{eq:inertia}) completes -- though does not improve~-- the well-known inequality
$$\int_{\R^d}p\,dx\le\frac{2I(0)}{dt^2}\,$$
which follows from the constancy (or the decay in the isentropic case) of
$$t\mapsto\int_{\R^d}\left(\frac12\,\rho|tu-x|^2+t^2\frac{dp}2\right)\,dx$$
for the mono-atomic gas.
\end{enumerate}

\subsection{Other ideal gases ($\gamma\ne\gamma_d$)}

As noted before, when $\gamma$ is not equal to $\gamma_d$, the transformation $(t,x,\rho,u,e)\mapsto(s,y,\bar\rho,v,\bar e)$ does not produce a gas flow, because the right-hand-side of (\ref{eq:tren}) does not vanish. Nevertheless $\bar S$ is divergence-free and we may apply Compensated Integrability to obtain (with Cauchy--Schwarz)
\begin{equation}
\label{eq:subCI}
\int_0^\tau ds\int_{\R^d}\bar\rho^\frac1d\bar p\,dy\leqslant c_dM^\frac1d\left(M\max(F_\alpha(0),F_\alpha(\tau))\right)^\frac12,
\end{equation}
where
$$F_\alpha(s)=\int_{\R^d}\left(\frac12\,\bar\rho|v|^2+\bar\rho\bar e\right)(s,y)\,dy.$$

If $\gamma<\gamma_d$, the integration of (\ref{eq:tren}) yields
$$\frac{dF_\alpha}{dt}\leqslant\frac{d\alpha}{(1-\alpha s)}\,(\gamma_d-\gamma)F_\alpha,$$
which gives
$$F_\alpha(s)=O((1-\alpha s)^{-K})=O((1+t\alpha)^K),\qquad K:=d(\gamma_g-\gamma).$$
Combining with (\ref{eq:subCI}), we conclude that

\begin{quote}
\begin{thm}\label{th:reel}
For an ideal gas with adiabatic constant $\gamma$ less than $\gamma_d$, the flows of finite mass, energy and moment of inertia satisfy
$$\int_0^Tt\,dt\int_{\R^d}\rho^\frac1dp\,dx=O\left(M^\frac1d\sqrt{M(E(0)+I(0))\,}\,T^{\frac{\gamma_d-\gamma}{\gamma_d-1}}\right).$$
\end{thm}
\end{quote}

\bigskip

Because of $\frac{\gamma_d-\gamma}{\gamma_d-1}<1$, the estimate above is still an improvement of (\ref{eq:estihp}) when $T\rightarrow+\infty$. For a gas whose molecules have $D$ freedom degrees ($D=5$ and $d=3$ for a di-atomic gas), this exponent is $1-\frac dD$\,.

Notice that the right-hand side is not homogeneous in terms of physical dimensions. It can be homogenized by rescaling, as in \cite{Ser_DPT}. The dimensional analysis suggests that the ultimate inequality be of the form
\begin{equation}\label{eq:ultime}
\int_0^Tt\,dt\int_{\R^d}\rho^\frac1dp\,dx\leqslant c_dM^\frac1d\sqrt{MI(0)\,}\,\left(T\sqrt{\frac{E(0)}{I(0)}\,}\right)^{\frac{\gamma_d-\gamma}{\gamma_d-1}}.
\end{equation}

\bigskip

When $\gamma>\gamma_d$ instead (not a realistic case), $s\mapsto F_\alpha(s)$ is non-increasing and we have simply
$$\int_0^\tau ds\int_{\R^d}\bar\rho^\frac1d\bar p\,dy\leqslant c_dM^\frac1d\sqrt{MF_\alpha(0)\,}\,.$$
This is exactly the same situation as in the mono-atomic case and we obtain again the estimate (\ref{eq:inertia}).

\section{Rarefied gases}\label{s:Bolt}

A mono-atomic rarefied gas is described by the particle distribution $fd\xi\,dx\,dt$, where $\xi$ is the velocity of particles. It obeys the Boltzman equation:
\begin{equation}
\label{eq:Boltz}
\partial_tf+\xi\cdot\nabla_xf=Q(f,f).
\end{equation}
The  bilinear operator $Q$ is given by
$$Q(f,f)(t,x,\xi)=\int_{\R^d}\!\int_{S^{d-1}}B(|\xi_1-\xi|,\omega)(f(t,x,\xi_1')f(t,x,\xi')-f(t,x,\xi_1)f(t,x,\xi))\,d\xi_1d\omega,$$
where as usual
$$\xi'=\xi+(\omega\cdot(\xi_1-\xi))\omega,\qquad\xi_1'=\xi-(\omega\cdot(\xi_1-\xi))\omega.$$
It is thus local in the variables $(t,x)$ but non-local in $\xi$.
The collision kernel $B$ is a non-negative function. The local conservation of mass, momentum and energy are encoded in the properties
$$\int_{\R^d}Q(g,g)\,d\xi=0,\qquad\int_{\R^d}Q(g,g)\xi\,d\xi=0,\qquad\int_{\R^d}Q(g,g)|\xi|^2d\xi=0,$$
for every $g=g(\xi)$ with reasonnable decay at infinity.

\bigskip

One defines classically the mass-momentum tensor
$$S(t,x)=\begin{pmatrix} \rho & m^T \\ m & T \end{pmatrix}:=\int_{\R^d}f(t,x,\xi)\binom1\xi\otimes\binom1\xi\,d\xi,$$
which is obviously symmetric, positive semi-definite. Integrating (\ref{eq:Boltz}) against $d\xi$ and $\xi\,d\xi$, we find formally that $S$ is Divergence-free. Integrating against $\frac12|\xi|^2d\xi$, we find also the conservation law
\begin{equation}
\label{eq:enB}
\partial_t\varepsilon+{\rm div}_x\vec q=0
\end{equation}
where
$$\varepsilon:=\int_{\R^d}f\,\frac{|\xi|^2}2\,d\xi,\qquad\vec q:=\int_{\R^d}f\,\frac{|\xi|^2}2\,\xi\,d\xi$$
are the energy density and the energy flux.

\bigskip

An existence theorem of {\em renormalized solutions} to the Cauchy problem has been established by R. DiPerna \& P.-L. Lions \cite{DPL} and completed by Lions \& N. Masmoudi \cite{LM2}. Under reasonnable assumptions on the collision kernel $B$, a renormalized solution exists whenever the initial data $f_0$ satisfies
\begin{equation}
\label{eq:ICB}
\int_{\R^d}\!\int_{\R_d}(1+|\xi|^2+|x|^2+|\log f_0|)f_0d\xi\, dx<+\infty.
\end{equation}
This assumption means that the total mass and energy, as well as the moment of inertia and the total entropy are finite at initial time.

In spite of the formal integration, one does not know (unless the space dimension equals one) whether $S$ is Divergence-free for renormalized solutions: according to \cite{LM2}, there exists indeed a $d\times d$  {\em defect tensor} $\Sigma$, which is symmetric positive semi-definite, such that the modified mass-momentum tensor
$$S_B=\begin{pmatrix} \rho & m^T \\ m & T+\Sigma \end{pmatrix}$$
is Divergence-free:
\begin{eqnarray*}
\partial_t\rho+{\rm div}_xm & = & 0, \\
\partial_t m+{\rm Div}_x(T+\Sigma) & = & 0.
\end{eqnarray*}
The first line implies, after an integration, that the total mass
$$M=\int_{\R^d}\!\int_{\R^d}f(t,x,\xi)\,d\xi\,dx$$
is a constant of the motion.
Regarding the conservation of energy, it is also known \cite{LM2} that the map
$$t\longmapsto\int_{\R^d}\!\int_{\R^d}f(t,x,\xi)\frac{|\xi|^2}2\,d\xi\,dx+\int_{\R^d}\frac12\Tr\Sigma\,dx$$
is constant, equal to the energy at initial time
$$E(0):=\int_{\R^d}\!\int_{\R^d}f_0(x,\xi)\frac{|\xi|^2}2\,d\xi\,dx.$$

\bigskip

Compensated Integrability and Cauchy--Schwarz give as usual
$$\int_0^{t_{\max}}dt\int_{\R^d}(\det S_B)^\frac1d\,dx\leqslant c_dM^\frac1d\sqrt{ME(0)\,}\,.$$
Because of $\Sigma\geqslant0_d$, which implies $0_{1+d}\leqslant S\leqslant S_B$ and thus $\det S\leqslant\det S_B$, the left-hand side dominates the integral of $(\det S)^\frac1d$.
\bigskip

We are now ready to apply our favorite projective transformation
$$s=\frac t{1+t\alpha}\,,\qquad y=\frac x{1+t\alpha}\,.$$
The new velocity variable is defined by
$$\chi=(1+t\alpha)\xi-\alpha x.$$
The new density is therefore
\begin{equation}
\label{eq:ffbar}
\bar f(s,y,\chi)=f\left(\frac s{1-s\alpha}\,,\frac y{1-s\alpha}\,,(1-s\alpha)\chi+\alpha y\right).\end{equation}
We observe that
\begin{eqnarray*}
(\partial_s+\chi\cdot\nabla_y)\bar f & = & (1+t\alpha)^2(\partial_t+\xi\cdot\nabla_x)f \\
& = & (1+t\alpha)^2Q(f,f)=:\bar Q_t(\bar f,\bar f).
\end{eqnarray*}

For most of kernels $B$, the function $\bar f$ is not a solution of the Boltzmann equation, because $\bar Q_t$ differs from $Q$. Actually, $\bar Q_t$ does depend explicitly upon the time variable, unless  $B(r,\omega)$ is of the form $r^{2-d}b(\omega)$. This exponent $2-d$ in the kinetic collision kernel corresponds to an inter-particles radial force with inverse power law $|x|^{-3}$, that is to the Calogero--Moser potential (see \cite{BV}).

What remains true in general is that $\bar Q_t(g,g)$ is annihilated by integrating against $\begin{pmatrix} 1 \\ \chi \\ |\chi|^2 \end{pmatrix}d\xi$, because it involves only the expression
$$\bar f(t,x,\chi_1')\bar f(t,x,\chi')-\bar f(t,x,\chi_1)\bar f(t,x,\chi)$$
for quadruplets $\chi_1',\chi',\chi_1,\chi$ compatible with the conservation of linear momentum and kinetic energy. Therefore $\bar f$ is expected to satisfy the same conservations and decay as $f$.

When taking the first moments of the density $\bar f$, we obtain the same quantities  as in Proposition \ref{p:projtr}~:
$$\bar\rho=\int_{\R^d}\bar f\,d\chi=(1+t\alpha)^d\rho,\qquad
\bar m=\int_{\R^d}\bar f\chi\,d\chi=(1+t\alpha)^d((1+t\alpha)m-\alpha\rho x)$$
and
$$\bar T=\int_{\R^d}\bar f\chi\otimes\chi\,d\chi=(1+t\alpha)^{d+2}T-\alpha(1+t\alpha)^{d+1}(m\otimes x+x\otimes m)+\alpha^2(1+t\alpha)^d\rho x\otimes x.$$
Proposition \ref{p:projtr} tells us that the positive tensor
$$\bar S_B=\begin{pmatrix} \bar\rho & \bar m^T \\ \bar m & \bar T+(1+t\alpha)^{d+2}\Sigma \end{pmatrix}$$
is Divergence-free in the coordinates $(s,y)$. Compensated Integrability thus yields
$$\int_0^{s_{\max}}ds\int_{\R^d}(\det\bar S)^\frac1ddy\leqslant\int_0^{s_{\max}}ds\int_{\R^d}(\det\bar S_B)^\frac1ddy\leqslant c_dM^\frac1d\sup_s\|\bar m(s,\cdot)\|_{\cal M}.$$
Since
$$\|\bar m(s,\cdot)\|_{\cal M}\leqslant\left(M\int_{\R^d}\bar f|\chi|^2d\chi\,dy\right)^\frac12,$$
and 
$$\int_{\R^d}\!\int_{\R^d}\bar f(s,y,\chi)\frac{|\chi|^2}2\,d\chi\,dy\,+\int_{\R^d}\frac12\Tr\bar\Sigma\,dy= E_\alpha(0):=\int_{\R^d}\!\int_{\R^d}\bar f_0(y,\chi)\frac{|\chi|^2}2\,d\chi\,dy,$$
we obtain
$$\int_0^{s_{\max}}ds\int_{\R^d}(\det\bar S)^\frac1ddy\leqslant c_dM^\frac1d\sqrt{ME_\alpha(0)\,}\,.$$
Going back to the original coordinates, where we have $\det \bar S=(1+t\alpha)^{d(d+3)}\det S$, this rewrites as
$$\int_0^{t_{\max}}(1+t\alpha)\,dt\int_{\R^d}(\det S)^\frac1ddx\leqslant c_dM^\frac1d\left(M\int_{\R^d}\!\int_{\R^d}f_0(x,\xi)\frac{|\xi-\alpha x|^2}2d\xi\,dx\right)^\frac12.$$
We conclude as in Paragraph \ref{p:mono} by letting $\alpha\rightarrow+\infty$, and then by optimizing the choice of the origin of $\R^n$. We thus obtain the following result.
\begin{thm}\label{th:Boltz}
Let the initial data $f_0\geqslant0$ satisfy the assumptions (\ref{eq:ICB}).
Then the renormalized solution satisfies
\begin{equation}
\label{eq:inBltz}
\int_0^{t_{\max}}t\,dt\int_{\R^d}(\det S)^\frac1ddx\leqslant c_dM^\frac1d\left(\frac14\,\int_{\R^d}\!\int_{\R^d}\rho_0(x)\rho_0(x')|x'-x|^2dx\,dx'\right)^\frac12.
\end{equation}
for an absolute constant $c_d$ that depends only upon the space dimension, and where $\det S$ is given in (\ref{eq:Andr}).
\end{thm}

\bigskip

 As noticed in \cite{Ser_DPT}, we may express
\begin{equation}
\label{eq:Andr}
\det S=\frac1{(d+1)!}\,\int_{\R^d}\cdots\int_{\R^d}f(t,x,\xi_0)\cdots f(t,x,\xi_d)V(\xi_0,\ldots,\xi_d)^2d\xi_0\cdots d\xi_d,
\end{equation}
where $V(\xi_0,\ldots,\xi_d)$ is the volume of the simplex spanned by the vertices $\xi_0,\ldots,\xi_d$ in $\R^d$.

\bigskip

Once again, (\ref{eq:inBltz}) seems to forbid source solutions of the Boltzman equation. If all the mass is concentrated at a single point at initial time, then the right-hand side vanishes, and this implies with (\ref{eq:Andr}) that for almost every $(t,x)$, 
$$a.e.\,\, (\xi_0,\ldots,\xi^d)\in(\R^d)^{1+d},\qquad (V(\xi_0,\ldots,\xi_d)>0)\Longrightarrow(f(t,x,\xi_0)\cdots f(t,x,\xi_d)=0).$$
In other words, the essential support of $f(t,x,\cdot)$ is contained in an affine hyperplane $\Pi_{t,x}$ of $\R^d$. We leave open the question whether such solutions exist besides $f\equiv0$. At least, they cannot belong to the class of renormalized solutions, since $f|\log f|$ is not integrable in space and velocity. The concentration of the support of $f(t,x,\cdot)$ can be interpreted as the fact that the temperature --~whatever this notion means for a flow off equilibrium~-- vanishes identically.

\subsection*{Appendix: Proof of (\protect\ref{eq:uncert})}

Since 
$$\int\!\int_{\R^{2d}}g(x)g(x')|x'-x|^2dx\,dx'=\min_{\hat x\in\R^d}2\int_{\R^d}g(x)\,dx\cdot\int_{\R^d}g(x)|x-\hat x|^2dx,$$
it is enough to prove
\begin{equation}
\label{eq:suffit}
\left(\int_{\R^d}g(x)\,dx\right)^{2+\frac2d}\leqslant c_d\int_{\R^d}g(x)|x|^2dx\cdot\int_{\R^d}g(x)^{1+\frac2d}dx
\end{equation}
Let $R>0$ be the radius of a ball $B_R$ centered at the origin, to be chosen later. We decompose the integral of $g$ as the sum of the integrals over $B_R$ and its complement. On the one hand, the H\"older inequality gives
$$\int_{B_R}g(x)\,dx\le \left(\int_{\R^d}g(x)^{1+\frac2d}dx\right)^{\frac d{d+2}}|B_R|^{\frac2{d+2}}= \left(\int_{\R^d}g(x)^{1+\frac2d}dx\right)^{\frac d{d+2}}(|B|R^d)^{\frac2{d+2}}.$$
On the other hand
$$\int_{B_R^c}g(x)\,dx\le \frac1{R^2}\,\int_{\R^d}g(x)|x|^2dx.$$
To balance both contributions, we choose
$$R^{4\frac{d+1}{d+2}}=\left(\int_{\R^d}g(x)^{1+\frac2d}dx\right)^{-\frac d{d+2}}\int_{\R^d}g(x)|x|^2dx$$
and we obtained the desired conclusion.

\begin{english}

\end{english}

\end{document}